\documentclass[preprint, fleqn, 11pt, 3p]{elsarticle}
\usepackage{color}
\usepackage{soul}
\usepackage[english]{babel}
\usepackage{graphicx}
\usepackage{amssymb,amsfonts,amsmath,calc,amstext}

\usepackage{mathrsfs}
\usepackage{amsmath}
\usepackage{graphicx}
\usepackage{epsfig}
\usepackage{subfig}
\usepackage{setspace}
\usepackage{tikz}
\usepackage{picture,xcolor}

\makeatletter
\newcommand{\rmnum}[1]{\romannumeral #1}
\newcommand{\Rmnum}[1]{\expandafter\@slowromancap\romannumeral #1@}
\makeatother

\newenvironment{proof}
 {\noindent{\bf Proof.}}

 \def\qed{\hbox{}\nobreak\hfill \vrule width 1.4mm height 1.4mm depth
0mm \par \goodbreak
\smallskip}

\newtheorem{theorem}{Theorem}

\newtheorem{proposition}{Proposition}
\newtheorem{remark}{Remark}

\newtheorem{corollary}{Corollary}

\begin{document}

\title{On the characterization of the structure of distributive uninorms}
\begin{frontmatter}

\author[jinan]{Wenwen Zong}
\ead{zongwen198811@163.com}

\author[suzhou]{Yong Su\corref{cor1}}
\ead{yongsu88@163.com}

\author[shandong]{Hua-Wen Liu}
\ead{hw.liu@sdu.edu.cn}

\tnotetext[22]{This work is supported by the National Natural Science Foundation of China (Grant nos.\,12271393, 12471435 and 12331016), the Natural Science Foundation of Shandong Province (Grant no.\,ZR2023MA092), Qing Lan Project of Jiangsu Province, and  Youth Innovative Research Team Plan of Higher Education of Shandong Province (Grant no.\,2023KJ103).}

\cortext[cor1]{Corresponding author.}

\address[jinan]{School of Mathematical Sciences, University of Jinan, Jinan, Shandong 250022, China}

\address[suzhou]{School of Mathematics Sciences, Suzhou University of Science and Technology, Suzhou, Jiangsu 215009, China}

\address[shandong]{School of Mathematics, Shandong University, 250100 Jinan, China}

\begin{abstract}
This paper focuses on distributive uninorms, which induce structures of commutative ordered semirings. We will show that the second uninorm must be locally internal on $A(e)$, and will present a complete characterization of the structure of such uninorms.
\end{abstract}


\end{frontmatter}

\section{Introduction, basic notions and results}

A \emph{uninorm} \cite{Yager96} is a binary operation $U: [0, 1]^2 \rightarrow [0, 1]$ such that $([0, 1], U, \leq)$ is a commutative, linearly ordered semigroup with neutral element $e \in [0, 1]$. A uninorm $U$ is known as a \emph{t-norm} whenever $e = 1$, while as a \emph{t-conorm} whenever $e = 0$. A  uninorm is called to be \emph{proper} if its neutral element belongs to $]0,1[$. Let $U$ be a proper uninorm. Then $U(0,1) \in\{0,1\}$, $U$ works as a t-norm in the square $[0,e]^2$, and as a t-conorm in $[e,1]^2$. Putting $A(e) = [0, 1]^2 \backslash ([0, e]^2 \cup [e, 1]^2)$, we then have $\min(x,y)\leq U(x,y) \leq \max(x,y)$ for any $(x, y) \in A(e)$. Denote by $U \equiv \langle T , e, S \rangle$ the class of all uninorms
with neutral element $e$, the underlying t-norm $T$ and the underlying t-conorm $S$.

In \cite{FodorYagerRybalov97}, the structure of uninorms was analyzed and three classes of uninorms were introduced: uninorms in $\mathcal{U}_{\min}$, uninorms in $\mathcal{U}_{\max}$, and representable uninorms. After then, several new classes of uninorms, i.e., uninorms continuous in the open unit square $]0, 1[^2$, idempotent uninorms, uninorms with continuous underlying operators and uninorms locally internal on $A(e)$, have appeared and been investigated. Here, we only recall uninorms locally internal on $A(e)$ and refer to  \cite{MasMassanetRuiz-AguileraTorrens15} for the other classes of uninorms. A uninorm $U$ is called to be \emph{locally internal on $A(e)$} if $U(x,y)\in\{x, y\}$ for any $(x, y) \in A(e)$. Denote by $\mathcal{U}_{\rm lin}$ the class of uninorms locally internal on $A(e)$.

Uninorms have proved to be useful in a wide range of fields like aggregation of information, expert systems, classifications, information measures, neural networks, pseudo-additive measures and approximate reasoning, and so on (see the recent survey \cite{MasMassanetRuiz-AguileraTorrens15}, monograph \cite{Drygas23} and the references therein for
more details). Due to a host of applications, uninorms have also been studied from the purely theoretical point of view. One of the research topics on uninorms is to characterize ones that verify certain properties. Such properties usually come from the concrete application. The distributivity is indispensable in nonlinear PDE \cite{BenvenutiMesiar04,Pap97} and several type of nonadditive integrals \cite[\S 5.6]{GrabischMarichalMesiarPap09}.  A uninorm $U_1$ is \emph{distributive} over a uninorm $U_2$ if, for all $x, y, z \in [0, 1]$,
\begin{eqnarray}\label{dis}
U_1(x, U_2(y,z)) = U_2(U_1(x, y), U_1(x, z)).
\end{eqnarray}
When a uninorm $U_1$ is distributive over a uninorm $U_2$, the algebra $([0, 1], U_2, U_1, \leq)$ is a commutative ordered semiring\footnote{A commutative semiring is an algebra $(S, +, \cdot)$ such that $(S, +)$ and $(S, \cdot)$ are commutative semigroups and $(a + b) c = ac + bc$ for all $a, b, c \in S$.  A commutative ordered semiring is a semiring $S$ equipped with a partial order $\leq$ such that both operations are increasing.}.

Research status of distributive uninorms is as follows:
\begin{itemize}
  \item The cases where both uninorms are in $\mathcal{U}_{\min}$ and $\mathcal{U}_{\max}$ were characterized in \cite{MasMayorTorrens02,MasMayorTorrens05}.
  \item  Ruiz and Torrens \cite{RuizTorrens03} characterized the structure of distributive idempotent uninorms.
  \item The distributivity of uninorms over disjunctive uninorms were discussed in \cite{RuizTorrens05,RuizTorrens09} under some partial continuity presupposition.
  \item The cases where one lies in anyone of the most well-known classes of uninorms and the other one is with rather weak assumptions were studied in \cite{SuLiuetal16,SuLiuetal18}.
\end{itemize}
For the studies above, a basic assumption is satisfied:
\begin{center}
\emph{at least one of uninorms lies in anyone of the most well-known classes of uninorms.}
\end{center}
This paper will characterize the structure of distributive uninorms without any assumption, going a step further than the above papers. For a couple of distributive uninorms, we will show that the second one must be locally internal on $A(e)$. The cases where the second uninorm is locally internal on $A(e)$ were discussed in \cite[Theorems 5, 10 and 11]{SuLiuetal18} under some additional internality, however, in this paper, we will give more concise presentations of the structure of distributive uninorms, without any assumption. On the contray, the decomposition method of distributive uninorms used in this paper is different from that used in \cite{SuLiuetal18}.

\section{Main results}
\label{sec:3}

In this section, the characterizations of the structure of distributive uninorms will be presented. Depending on the order relationships of the neutral elements $e_1, e_2$ of two uninorms $U_1, U_2$, we will distinguish three cases: $e_1 = e_2$, $e_1 > e_2$ and $e_1 < e_2$.

\subsection{The case $e_1 = e_2$}

This case has been completely characterized in \cite[Theorem 5]{SuLiuetal18}. Here, we present it for the sake of completeness.

\begin{theorem} {\rm \cite[Theorem 5]{SuLiuetal18}}
Let $U_1$ and $U_2$ be two uninorms with the same neutral element $e$. Then
$U_1$ is distributive over $U_2$ if and only if $U_2$ is idempotent\footnote{A uninorm $U$ is idempotent if $U(x, x) = x$ for all $x \in [0, 1]$.} and
$U_1(x,y) = U_2(x,y) \in \{x, y\}$ for all $(x,y) \in A(e)$.
\end{theorem}

\subsection{The case $e_1 > e_2$}

\begin{proposition}\label{mainpro}
Let $0 < e_2 < e_1 < 1$ and let $U_1\equiv\langle T_1,e_1,S_1\rangle$ be distributive over $U_2\equiv\langle T_2,e_2,S_2\rangle$. Then
\begin{itemize}
  \item[{\rm ({\rmnum 1})}] $T_2 = \min$.
  \item[{\rm ({\rmnum 2})}] $U_1(x, y) = \min(x, y)$ for any $x\in [0, e_2]$ and $y\in [e_2, e_1]$.
  \item[{\rm ({\rmnum 3})}] $U_2(x,y) = \min(x, y)$ for any $x \in \,[0, e_2[$ and $y \in [e_2, e_1]$.
  \item[{\rm ({\rmnum 4})}] $U_1(x,y) = U_2(x,y)\in \{x, y\}$ for any $x \in \,[0, e_2[$ and $y \in [e_1, 1]$. In particular, $U_2(y_0,y_0)=y_0$ whenever $U_2(x_0,y_0)=y_0$ for some $x_0 \in [0, e_2[$ and $y_0 \in [e_1, 1]$.
  \item[{\rm ({\rmnum 5})}] there exists a uninorm $\widetilde{U}$ with neutral element $\frac{e_1-e_2}{1-e_2}$ such that, for any $x, y \in [e_2, 1]$,
  \begin{eqnarray}\label{U1=U}
  U_1(x,y) =  e_2 + (1 - e_2) \widetilde{U}\left (\frac{x - e_2}{1 - e_2},  \frac{y - e_2}{1 - e_2}\right),
  \end{eqnarray}
 and $\widetilde{U}$ is distributive over~$S_2$.
\end{itemize}
\end{proposition}

\begin{proof}
({\rmnum 1}) Considering arbitrary fixed $x\in [0, e_2]$ and $y\in [e_2, e_1]$. Then $U_1(x,e_2) \leq x$ and the distributivity yields that
\begin{eqnarray}\label{x=U_2(U_1(x,e_2),x)}
x = U_1(x,U_2(e_2, e_1)) =U_2(U_1(x,e_2),U_1(x,e_1)) = U_2(U_1(x,e_2),x) \leq U_2(x, x) \leq x,
\end{eqnarray}
implying $U_2(x, x) =x$, or equivalently $T_2=\min$.

({\rmnum 2}) From Eq.\,\eqref{x=U_2(U_1(x,e_2),x)} and $T_2=\min$, it follows that $x = U_2(U_1(x,e_2),x) = U_1(x, e_2)$ for all $x \in [0, e_2]$. Since $U_1$ is increasing and has neutral element $e_1$, we obtain that $U_1(x, y) = \min(x, y)$ for any $x\in [0, e_2]$ and $y\in [e_2, e_1]$.

({\rmnum 3}) Suppose $x \in \,[0, e_2[$ and $y \in [e_2,e_1]$. We then have $x \leq U_2(x,y) \leq y \leq e_1$. Item ({\rmnum 2}) yields $U_1(e_2,x) = x$ and $U_1(e_2,y)= e_2$. Hence,
\begin{eqnarray*}
U_1(e_2,U_2(x, y)) = U_2(U_1(e_2, x), U_1(e_2, y)) = x,
\end{eqnarray*}
implying $U_2(x,y) < e_2$ (otherwise, $x = U_1(e_2,U_2(x, y)) = e_2$ by ({\rmnum 2})). Consequently, ({\rmnum 2})) yields  $x = U_1(e_2,U_2(x, y)) = U_2(x,y)$.

({\rmnum 4}) Suppose $x \in [0, e_2[$ and $y \in [e_1, 1]$. Then $x \leq U_2(x,y) \leq y$. The distributivity and ({\rmnum 2}) yield
\begin{align}\label{2}
U_1(x,y) = U_1(y,U_2(x, e_1)) = U_2(U_1(y, x), U_1(y, e_1)) = U_2((U_1(y, x), y). \\\label{3}
 U_1(x,y) = U_1(x,U_2(y, e_2)) = U_2(U_1(x, y), U_1(x, e_2)) = U_2((U_1(x, y), x)
\end{align}
\begin{itemize}
  \item If $U_1(x,y)\leq e_2$,  then Eq.\eqref{3} indicates $x \leq U_1(x,y) =  U_2((U_1(x, y), x) \leq x,$ implying $U_1(x,y) = x$. By Eq.\eqref{2}, we have $U_2(x,y) = x$.
  \item If $U_1(x,y) > e_2$, then Eq\,.\eqref{2} yields $y \leq U_2(U_1(y, x), y) = U_1(y,x) \leq y,$ i.e., $U_1(x,y) = y$ and $U_2(y,y)=y$, which, together with  Eq.\eqref{3}, gives $U_2(x,y) = y$.
\end{itemize}

({\rmnum 5}) Item ({\rmnum 2}) says $e_2$ is an idempotent element of $U_1$. Hence, the following binary function $\widetilde{U}: [0, 1]^2 \rightarrow [0, 1]$ defined by
\begin{eqnarray*}
\widetilde{U}(x, y) = \frac{U_1(e_2 + (1 - e_2)x, e_2 + (1 - e_2)y) - e_2}{1 - e_2}
\end{eqnarray*}
is a uninorm $\widetilde{U}$ with neutral element $\frac{e_1-e_2}{1-e_2}$ and Eq.\,\eqref{U1=U} holds.  Consider $x, y, z \in [e_2, 1]$,  the distributivity yields
\begin{eqnarray*}
S_2 \Bigg(\widetilde{U}\Big(\frac{x - e_2}{1 - e_2}, \frac{y - e_2}{1 - e_2}\Big), \widetilde{U}\Big(\frac{x - e_2}{1 - e_2}, \frac{z - e_2}{1 - e_2}\Big)\Bigg) = \widetilde{U}\Bigg(\frac{x - e_2}{1 - e_2}, S_2\Big(\frac{y - e_2}{1 - e_2}, \frac{z - e_2}{1 - e_2}\Big)\Bigg),
\end{eqnarray*}
i.e., $\widetilde{U}$ is distributive over~$S_2$.\qed
\end{proof}

From the previous proposition, we can deduce the following result:
\begin{corollary}
Let $0 < e_2 < e_1 < 1$. If $U_1\equiv\langle T_1,e_1,S_1\rangle$ is distributive over $U_2\equiv\langle T_2,e_2,S_2\rangle$, then $U_2 \equiv\langle \min, e_2, S_2\rangle \in \mathcal{U}_{\rm lin}$.
\end{corollary}

\begin{theorem}\label{theorem1}
Let $0 < e_2 < e_1 < 1$, $U_1\equiv\langle T_1,e_1,S_1\rangle$ and $U_2\equiv\langle \min,e_2,S_2\rangle \in {\mathcal{U}_{\rm lin}}$. Then $U_1$ is distributive over $U_2$ if and only if
\begin{itemize}
  \item[{\rm ({\rmnum 1})}] $U_1(x, y) = U_2(x, y) \in \{x, y\}$ for any $x \in [0, e_2[$ and $y \in [e_2, 1]$. In particular, $U_2(y_0,y_0)=y_0$ when $U_2(x_0,y_0)=y_0$ for some $x_0 \in [0, e_2[$ and $y_0 \in [e_1, 1]$.
   \item[{\rm ({\rmnum 2})}] $U_1(x, y) = U_2(x, y) = \min(x, y)$ for any $x \in [0, e_2[$ and $y \in [e_2, e_1]$.
  \item[{\rm ({\rmnum 3})}] there exists a uninorm $\widetilde{U}$ with neutral element $\frac{e_1-e_2}{1-e_2}$ such that, for any $x, y \in [e_2, 1]$,
  \begin{eqnarray*}
  U_1(x,y) =  e_2 + (1 - e_2) \widetilde{U}\left (\frac{x - e_2}{1 - e_2},  \frac{y - e_2}{1 - e_2}\right),
  \end{eqnarray*}
 and $\widetilde{U}$ is distributive over~$S_2$.
\end{itemize}
\end{theorem}

\begin{proof}
The result necessarily follows from Proposition \ref{mainpro} and the converse can be proved by verifying case-by-case, which is similar to that of Theorem 5 in \cite{SuLiuetal18}.  \qed
\end{proof}

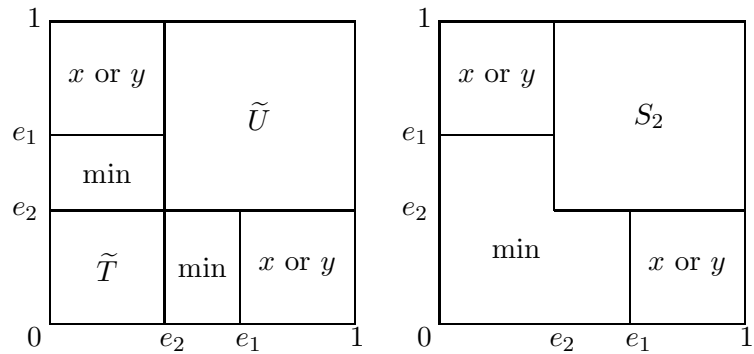
\begin{figure}[ht]
\begin{center}
\setlength{\unitlength}{1.0 cm}
\begin{picture}(5,5)(0,0)
  \put(0.5,0.5){\framebox(4,4){}}
  \put(0.5,0.5){\makebox(1.5,1.5){$\widetilde{T}$}}
  \put(0.5,2){\makebox(1.5,1){$\min$}}
   \put(2,0.5){\makebox(1,1.5){$\min$}}
\put(3,0.5){\makebox(1.5,1.5){$x$ or $y$}}
\put(0.5,3){\makebox(1.5,1.5){$x$ or $y$}}

  \put(2,2){\makebox(2.5,2.5){$\widetilde{U}$}}

  \put(0.5,2){\line(1,0){4}}
  \put(2,0.5){\line(0,1){4}}
\put(0.5,3){\line(1,0){1.5}}
  \put(3,0.5){\line(0,1){1.5}}

  \put(0.2,0.2){$0$}
\put(1.94,0.2){$e_2$}
\put(2.94,0.2){$e_1$}
\put(4.44,0.2){$1$}
\put(0,2.94){$e_1$}
\put(0,1.94){$e_2$}
\put(0.2,4.44){$1$}
 \end{picture}
\begin{picture}(5,5)(0,0)
  \put(0.5,0.5){\framebox(4,4){}}
  \put(0.5,0.5){\makebox(2,2){$\min$}}
  \put(2,2){\makebox(2.5,2.5){${S_2}$}}
   \put(2,2){\line(1,0){2.5}}
   \put(2,2){\line(0,1){2.5}}
  \put(0.5,3){\line(1,0){1.5}}
  \put(3,0.5){\line(0,1){1.5}}
\put(3,0.5){\makebox(1.5,1.5){$x$ or $y$}}
\put(0.5,3){\makebox(1.5,1.5){$x$ or $y$}}

  \put(0.2,0.2){$0$}
\put(1.94,0.2){$e_2$}
\put(2.94,0.2){$e_1$}
\put(4.44,0.2){$1$}
\put(0,2.94){$e_1$}
\put(0,1.94){$e_2$}
\put(0.2,4.44){$1$}
 \end{picture}
 \caption{The structure of uninorms $U_1$ and $U_2$ when $e_2<e_1$.}\label{Fig.e_2-e_11}
 \end{center}
\end{figure}

For the case that $0 < e_2 < e_1 < 1$, the general structure of $U_1$ and $U_2$ from Theorem~\ref{theorem1} can be viewed in Fig.~\ref{Fig.e_2-e_11}.

\begin{remark}
Let $U_1\equiv\langle T_1,e_1,S_1\rangle$ be distributive over $U_2\equiv\langle T_2,e_2,S_2\rangle \in {\mathcal{U}_{\rm lin}}$. Then the second uninorm $U_2$ has a characterizing operator $g_2$. \cite[Theorem 11]{SuLiuetal18} presented such a couple of uninorms depending on the order relationships between $a$ and $e_1$, where $a=\sup\{z\mid g_2(z)=e_2\}$:
\begin{itemize}
  \item The case $a = e_1$ leads that $U_2$ is idempotent and the characterization of such uninorms was presented in \cite[Theorem 8]{SuLiuetal18}.
  \item For the case $a < e_1$, under the assumption that $U_1(x,a)=x$ for all $x\geq a$, such uninorms was characterized.
\end{itemize}
In our work, we provide a comprehensive characterization of the structure of distributive uninorms without any assumptions. Moreover, it is easy to see that the decomposition method of distributive uninorms used in this paper is different from that used in \cite[Theorem 11]{SuLiuetal18}.
\end{remark}

\subsection{The case $e_1 < e_2$}

Here, we only list the corresponding results avoiding their proofs because one can prove them by a procedure completely similar to the case  $e_1 > e_2$.

\begin{proposition}
Let $0 < e_1 < e_2 < 1$ and let $U_1\equiv\langle T_1,e_1,S_1\rangle$ be distributive over $U_2\equiv\langle T_2,e_2,S_2\rangle$. Then
\begin{itemize}
  \item[{\rm ({\rmnum 1})}] $S_2 = \max$.
  \item[{\rm ({\rmnum 2})}] $U_1(x, y) = \max(x, y)$ for any $x \in \,[e_2,1]$ and $y \in [e_1, e_2]$.
  \item[{\rm ({\rmnum 3})}] $U_2(x, y) = \max(x, y)$ for any $x \in \,]e_2,1]$ and $y \in [e_1, e_2]$.
  \item[{\rm ({\rmnum 4})}] $U_1(x, y) = U_2(x, y) \in \{x, y\}$ for any $x \in ]e_2,1]$ and $y \in [0, e_1]$. In particular, $U_2(y_0,y_0)=y_0$ when $U_2(x_0,y_0)=y_0$ for some $x_0 \in ]e_2,1]$ and $y_0 \in [0,e_1]$.
  \item[{\rm ({\rmnum 5})}] there exists a uninorm $\widehat{U}$ with neutral element $\frac{e_1}{e_2}$ such that
  \begin{eqnarray*}U_1(x,y) =  e_2 \widehat{U}\left (\frac{x }{ e_2},  \frac{y }{ e_2}\right) \end{eqnarray*}
   for any $x, y \in [0,e_2]$, and $\widehat{U}$ is distributive over~$T_2$.
\end{itemize}

\end{proposition}

\begin{corollary}
Let $0 < e_1 < e_2 < 1$ and $U_1\equiv\langle T_1,e_1,S_1\rangle$ be distributive over $U_2\equiv\langle T_2,e_2,S_2\rangle$. Then $U_2\equiv\langle T_1, e_2, \max \rangle\in \mathcal{U}_{\rm lin}$.
\end{corollary}

\begin{theorem}\label{theorem2}
Let $0 < e_1 < e_2 < 1$, $U_1\equiv\langle T_1,e_1,S_1\rangle$ and $U_2\equiv\langle T_2,e_2, \max \rangle \in {\mathcal{U}_{\rm lin}}$. Then $U_1$ is distributive over $U_2$ if and only if the following statements hold:
\begin{itemize}
  \item[{\rm ({\rmnum 1})}] $U_1(x, y) = U_2(x, y) \in \{x, y\}$ for any $x \in ]e_2,1]$ and $y \in [0, e_2]$. In particular, $U_2(y_0,y_0)=y_0$ when $U_2(x_0,y_0)=y_0$ for some $x_0 \in ]e_2,1]$ and $y_0 \in [0,e_1]$.
  \item[{\rm ({\rmnum 2})}] $U_1(x, y) = U_2(x, y) = \max(x, y)$ for any $x \in \,]e_2,1]$ and $y \in [e_1, e_2]$.
  \item[{\rm ({\rmnum 3})}] there exists a uninorm $\widehat{U}$ with neutral element $\frac{e_1}{e_2}$ such that
  \begin{eqnarray*}U_1(x,y) =  e_2 \widehat{U}\left (\frac{x }{ e_2},  \frac{y }{ e_2}\right) \end{eqnarray*}
   for any $x, y \in [0,e_2]$, and $\widehat{U}$ is distributive over~$T_2$.
  \end{itemize}
\end{theorem}

\begin{figure}[ht]
\begin{center}
\setlength{\unitlength}{1.0 cm}
\begin{picture}(5,5)(0,0)
  \put(0.5,0.5){\framebox(4,4){}}
  \put(3,3){\makebox(1.5,1.5){$\widehat{S}$}}
  \put(0.5,3){\makebox(1.5,1.5){$x$ or $y$}}
   \put(3,0.5){\makebox(1.5,1.5){$x$ or $y$}}
\put(3,2){\makebox(1.5,1){$\max$}}
\put(2,3){\makebox(1,1.5){$\max$}}
\put(2,3){\line(0,1){1.5}}
  \put(3,2){\line(1,0){1.5}}

  \put(0.5,0.5){\makebox(2.5,2.5){$\widehat{U}$}}

  \put(0.5,3){\line(1,0){4}}
  \put(3,0.5){\line(0,1){4}}

  \put(0.2,0.2){$0$}
\put(1.94,0.2){$e_1$}
\put(2.94,0.2){$e_2$}
\put(4.44,0.2){$1$}
\put(0,2.94){$e_2$}
\put(0,1.94){$e_1$}
\put(0.2,4.44){$1$}
 \end{picture}
\begin{picture}(5,5)(0,0)
  \put(0.5,0.5){\framebox(4,4){}}
  \put(0.5,0.5){\makebox(2.5,2.5){$T_2$}}
  \put(3,3){\makebox(1.5,1.5){${\max}$}}
   \put(0.5,3){\line(1,0){2.5}}
   \put(3,0.5){\line(0,1){2.5}}
\put(0.5,3){\makebox(1.5,1.5){$x$ or $y$}}
   \put(3,0.5){\makebox(1.5,1.5){$x$ or $y$}}

\put(2,3){\line(0,1){1.5}}
  \put(3,2){\line(1,0){1.5}}

  \put(0.2,0.2){$0$}
\put(1.94,0.2){$e_1$}
\put(2.94,0.2){$e_2$}
\put(4.44,0.2){$1$}
\put(0,2.94){$e_2$}
\put(0,1.94){$e_1$}
\put(0.2,4.44){$1$}
 \end{picture}
 \caption{The structure of uninorms $U_1$ and $U_2$ when $e_1<e_2$.}\label{Fig.e_1-e_21}
 \end{center}
\end{figure}
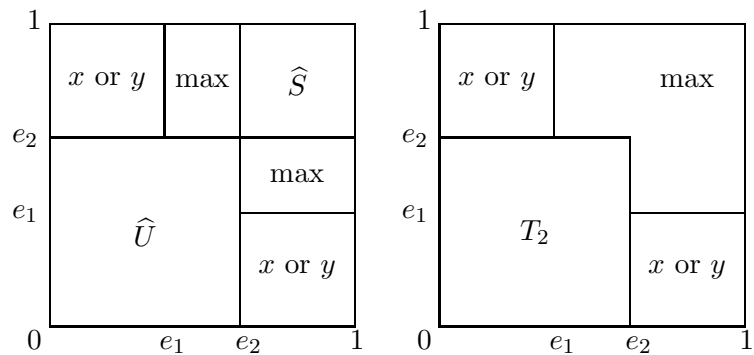

For the case that $0 < e_1 < e_2 < 1$, the general structure of $U_1$ and $U_2$ from Theorem~\ref{theorem2} can be viewed in Fig.~\ref{Fig.e_1-e_21}.

\bibliographystyle{plain}

\begin{thebibliography}{99}
%
\bibitem{BenvenutiMesiar04}
P. Benvenuti, R. Mesiar, Pseudo-arithmetical operations as a basis for the general measure and integration theory, Inf. Sci.,  160 (2004) 1-11.

%
\bibitem{Drygas23}
P. Dryga\'{s}, Uninorms and their Applications. Wydawnictwo Uniwersytetu Rzeszowskiego, Rzesz\'{o}w, 2023.
%
\bibitem{FodorYagerRybalov97}
J. Fodor, R.R. Yager, A. Rybalov, Structure of uninorms, Int. J. Uncertain. Fuzziness Knowl.-Based Syst. 5(4) (1997) 411--427.
%
\bibitem{GrabischMarichalMesiarPap09}
M. Grabisch, J.-L. Marichal, R. Mesiar, E. Pap, Aggregation Functions, Cambridge University Press, 2009.

\bibitem{MasMayorTorrens02}
M. Mas, G. Mayor, J. Torrens, The distributivity condition for uninorms and t-operators, Fuzzy Sets Syst. 128(2) (2002) 209-225.
%
\bibitem{MasMayorTorrens05}
M. Mas, G. Mayor, J. Torrens, Corrigendum to ``the distributivity condition for uninorms and t-operators", Fuzzy Sets Syst. 153 (2005) 297-299.
%
\bibitem{MasMassanetRuiz-AguileraTorrens15}
M. Mas, S. Massanet, D. Ruiz-Aguilera, J. Torrens, A survey on the existing classes of uninorms, J. Intell. Fuzzy Syst. 29 (2015) 1021-1037.


%
\bibitem{Pap97}
E. Pap, Decomposable measures and nonlinear equations, Fuzzy Sets Syst., 92 (1997) 205-221.
%
\bibitem{RuizTorrens03}
D. Ruiz, J. Torrens, Distributive idempotent uninorms, Int. J. Uncertain. Fuzziness Knowl.-Based Syst. 11(4) (2003) 413-428.
%
\bibitem{RuizTorrens05}
D. Ruiz-Aguilera, J. Torrens, Distributivity of strong implications over conjunctive and disjunctive uninorms, Kybernetika 42 (2005) 319-336.
%
\bibitem{RuizTorrens09}
D. Ruiz-Aguilera, J. Torrens, S- and R-implications from uninorms continuous in $]0, 1[^2$ and their distributivity over uninorms, Fuzzy Sets Syst. 160 (2009) 832-852.
%
\bibitem{SuLiuetal16}
Y. Su, H. Liu, D. Ruiz-Aguilera, J. V. Riera, J. Torrens, On the distributivity property for uninorms, Fuzzy Sets Syst. 287 (2016) 184-202.

%
\bibitem{SuLiuetal18}
Y. Su, H. Liu, J. V. Riera, D. Ruiz-Aguilera,  J. Torrens, The distributivity equation for uninorms revisited, Fuzzy Sets Syst. 334 (2018) 1-23.
%
\bibitem{Yager96}
R. R. Yager, A. Rybalov, Uninorm aggregation operators, Fuzzy Sets Syst. 80 (1996) 111-120.

\end{thebibliography}

\end{document}